\documentclass[cmp,numbook]{svjour}

\usepackage{amsmath}
\usepackage{amsfonts,amssymb}
\usepackage{psfig}

\journalname{Communications in Mathematical Physics}
\newtheorem{Theorem}{Theorem}

\newtheorem{Lemma}[Theorem]{Lemma}

\newcommand{\R}{\mathbf{R}}
\newcommand{\Z}{\mathbf{Z}}
\newcommand{\N}{\mathbf{N}}
\newcommand{\I}{\mathbf{I}}
\newcommand{\U}{\mathbf{U}}
\newcommand{\tr}{\mbox{tr}}

\begin{document}

\title{Multi-scaling of the $n$-point density function for coalescing Brownian motions}
\author{Ranjiva Munasinghe\inst{1}\thanks{Electronic address: ranm@maths.warwick.ac.uk}\and R. Rajesh\inst{2}
\thanks{Electronic address: rrajesh@imsc.res.in}
\and Roger Tribe\inst{1}\thanks{Electronic address: tribe@maths.warwick.ac.uk}
\and Oleg Zaboronski\inst{1}\thanks{Electronic address: olegz@maths.warwick.ac.uk}}

 \institute{Mathematics Institute, University of Warwick, Gibbet Hill Road, Coventry CV4 7AL, UK
\and  
 Institute of Mathematical Sciences, CIT Campus, Taramani, 
Chennai-600113, India}

\date{\today}
\maketitle

\begin{abstract}
This paper gives a derivation for the large time asymptotics of the $n$-point density function
of a system of coalescing Brownian motions on $\R$.
\end{abstract}

\section{Introduction and statement of the main result} \label{s1}
The single species reaction-diffusion systems $A+A \rightarrow A$ (coalescence) and
$A+A \rightarrow 0$ (annihilation) have been studied extensively in
recent times \cite{Cardy,mass1,mass2,ben,Doer,Bram,Howd,Lee}.
A common interest in these papers is the departure from mean field  behaviour
and the computation of exact long-term asymptotics for the particle density.

Recently, it was predicted in  \cite{Mun1} that the large time
asymptotics for the probability $P_t(n,\Delta V)$ of finding $N$ particles in
a fixed volume $\Delta V$:
\begin{equation}
\label{eq:mainpaper1}
P_t (n, \Delta V) \sim 	
\left\{ \begin{array}{ll}
t^{-\frac{n}{2}-\frac{n(n-1)}{4}} & d = 1 \\
\left(\frac{\ln t}{t}\right)^n (\ln t)^{-\frac{n(n-1)}{2}} & d = 2 \\
t^{-n} & d > 2 
\end{array} \right.		
\end{equation}  
Note that the predictions for $d>2$ agree with mean field behaviour.
The second part of the exponent in $d = 2$ reflects multi-scaling, or deviation from
linear scaling. In $d=2$ the multi-scaling is manifested in the second logarithmic term. 
This type of scaling is indicative of particles being {\it anti-correlated} \cite{Mun1}.
These predictions were obtained by use of the dynamical perturbative renormalization group
methods in a field theoretic setting \cite{Cardy}.
The setting here is for finite rate reactions, leading to annihilations of randomly walking
particles on a fixed lattice. After renormalization the
large time limit rate tends to an instantaneous reaction \cite{Cardy,Lee}.
Moreover the predictions carry over to the coalescing case, since they have the
same effective field theory \cite{HH,mass2}.

The aim of this paper is to verify the conjectures of
\cite{Mun1} in dimension $d=1$, where detailed probabilistic tools are available.
We consider a system of coalescing Brownian particles on the real line. 
Each particle evolves independently until it collides with another particle, at which time 
the two colliding particles instantaneously coalesce into one. The $n$-point density function 
is defined, for distinct $y_1,\ldots,y_n$, by
\[
P \left[ \mbox{there exist particles in $dy_1,\ldots,dy_n$ at time $t$} \right]
= \rho_n(y_1,\ldots,y_n;t) \, dy_1 \ldots dy_n.
\]
Note that $\rho_n$ depends on the initial particle distribution. The existence
of such a density is dicussed in the appendix. 

Our main result can be stated precisely as the following asymptotic:
\begin{equation} \label{eq:mainconj}
\rho_n(y_1,\ldots,y_n;t) \sim t^{-\alpha(n)} \prod_{1 \leq i < j \leq n} |y_i-y_j|
\end{equation}where $\alpha(n) = n/2 + n(n-1)/4$. 
This asymptotic has the meaning that the left hand side is
bounded above and below by constant multiples $C_1,C_2$ of the right hand side.
The upper bound is not actually an asymptotic, in that it 
holds simultaneously for all initial conditions, for 
all $t >0$ and for all $|y_i| \leq L t^{1/2}$, 
with a constant $C_2(L,n)< \infty$ depending only on $n$ and $L$.
The lower bound for all $t \geq t_0$ and all $|y_i| \leq L t^{1/2}$, with 
a constant $C_1(n,L,t_0)>0$ that depends on $n$,$L$,$t_0$ and also 
on the initial condition. For the lower bound, the  
initial condition must satisfy a mild non-degeneracy condition; in particular it 
holds for deterministic initial conditions provided the gap between successive particles is bounded, 
but also if we assume the set of initial positions of particles 
$\{X^i_0:i \geq 1\}$ is non-zero, translationally invariant and spatially ergodic, 
in the sense that the distribution of $\{X^i_0:i \geq 1\} \cap (-\infty,0]$ and 
$\{X^i_0:i \geq 1\} \cap [L,\infty)$ become independent as $L \to \infty$. 

The key tool is the Karlin-McGregor formula for the non-coincidence 
probabilities for Brownian motions. Useful upper and lower bounds
on this transition density, which already display the
key anomalous scaling term $t^{-\alpha(n)}$, are developed in section \ref{s3}, by exploiting
a representation known as the Harish-Chandra-Itzykson-Zuber
formula (developed for random matrix problems). 

The empty interval method, and its generalizations, have been used to derive
expressions for higher order correlation functions in \cite{mass1,mass2,ben}
for one-dimensional systems with instantaneous reactions.
Large time asymptotics for the $n$-point
correlation function density for the coalescing case with
Poissonian initial conditions are found in \cite{ben}, while the corresponding
results for the annihilating case are given in \cite{mass2}.
It was shown in this special case that the $n$-point density correlation function for the two systems
are the same apart from the amplitude. This set of exact results was used to test the
predictions (\ref{eq:mainpaper1}) in \cite{Mun1}. The large time scaling of the formulae for the
$n$-point density functions given in \cite{mass2,ben} are not obvious
as they involve a large combinatorial sum of terms with alternating signs.

The problem of deriving rigourously the logarithmic
corrections (\ref{eq:mainpaper1}) to the mean field answers in dimension $d=2$ remains open. 
It would also be interesting to find out if there is a natural multi-fractal
interpretation of the multi-scaling. 
Another simple system for which RG calculations predict multi-scaling in the stationary state 
is the system of aggregating massive point clusters with stationary source
of light particles. This system is relevant to turbulence, see \cite{Oleg1}. 
It would be interesting to generalize the methods of present paper
to prove multi-scaling for such cluster-cluster aggregation.
\section{Proof of the main result} \label{s2}
The proof is based on the following two lemmas. 
The first is a bound on the Karlin-McGregor formula for the  
transition density for non-intersecting Brownian motions \cite{Karlin,Mohanty}.
Fix $x=(x_1,\ldots,x_n) \in \R^n$ and let 
$(X^{x_i}: i=1,\ldots,n)$ be independent Brownian motions with $X^{x_i}_0=x_i$. 
For $y \in \R^n$, denote by $G^{KM}_t(x,y)$ the density of the probability measure 
\[ 
P \left[ X^{x_i}_t \in dy_i \; \mbox{and the paths $(X^{x_i}_s:s \in [0,t])$ are 
non-intersecting for $i=1,\ldots,n$.} \right].
\]
Then the Karlin-McGregor formula is 
\begin{eqnarray} \label{KM}
G_t^{KM}(x,y)
& = & \left| \begin{array}{ccc}
G_t(x_1,y_1) & \ldots & G_t(x_1,y_n) \\
\vdots && \vdots\\
G_t(x_n,y_1) & \ldots & G_t(x_n,y_n)
\end{array}
\right|
\end{eqnarray}
where $G_t(a,b) = (2 \pi t)^{-1/2} e^{-|a-b|^2/2t}$ is the one dimensional 
Brownian transition density. 
Note that the function $G^{KM}_t(x,y)$ is the transition density for an $n$-dimensional
Brownian motion killed on the set $\cup_{i \neq j}\{y_i=y_j\}$. One can check directly
that the determinant expression (\ref{KM}) satisfies the heat equation with zero
Dirichlet boundary conditions on this boundary, and that $\lim_{t \downarrow 0}
G^{KM}_t(x,y) = \delta_{x=y}$.

The following lemma, proved in section \ref{s3}, gives usable bounds on $G^{KM}_t(x,y)$.
\begin{Lemma} \label{l1}
For all $x_1< \ldots <x_n$, $y_1< \ldots<y_n$ and $t>0$ 
\begin{equation} \label{KMbounds}
c_n \prod_{i=1}^{n} G_t(x_i,y_{n-i+1})
\leq \frac{G_t^{KM}(x,y)}{\Delta(xt^{-1/2}) \, \Delta(yt^{-1/2})} 
 \leq c_n \prod_{i=1}^{n} G_t (x_i,y_i),
\end{equation}
where $c_n^{-1} = \prod_{i=1}^n i!$ and $\Delta(x)$ is the 
Vandermonde determinant defined by 
\[
\Delta(x) = \Delta(x_1,\ldots,x_n) = \prod_{1\leq i<j\leq n} (x_i-x_j).
\] 
\end{Lemma}

The second lemma, proved in section \ref{s4}, gives a simple upper 
bound on the $n$-point correlation function, which reflects the 
intuition that particles should be anti-correlated, in
that the presence of a particle in $dx$ decreases the likelihood 
that there is a particle at another point $dy$. 
\begin{Lemma} \label{l2}
For any initial distribution of particles, the $n$-point density function satisfies
\[
\rho_n(y_1,\ldots,y_n;t) \leq (\pi t)^{-n/2} \quad \mbox{for all $y_i \in \R$ and $t>0$.} 
\]
\end{Lemma}
 
The lower bound in the asymptotic (\ref{eq:mainconj}) follows quickly from the
lower bound on the Karlin-McGregor formula (\ref{KMbounds}). Indeed, list the 
set of initial positions of particles as $\{X^i_0:i \geq 1\}$ and let
$\Omega_0(t)$ be the event that there exist particles $X^{i_1}_0, \ldots,X^{i_n}_0$  
satisfying 
\begin{equation} \label{lbparticles}
X^{i_k}_0 \in [2kt^{1/2},(2k+1)t^{1/2}] \quad \mbox{for $k=1,\ldots,n$.}
\end{equation} 
Then, applying the Karlin-McGregor formula between
the points $X^{i_1}_0, \ldots, X^{i_n}_0$ and $y_1, \ldots,y_n$ we 
find, for $t \geq 1$, 
\begin{eqnarray*}
&& \hspace{-.4in} 
P \left[ \mbox{there exist particles at $dy_1,\ldots,dy_n$ at time $t$} \right] \\
& \geq & E \left[ \I(\Omega_0(t)) \, 
G^{KM}_t\left( (X^{i_1}_0,\ldots,X^{i_n}_0), (y_1,\ldots,y_n) \right) \right] \\
& \geq & c_n E \left[ \I(\Omega_0(t)) 
\Delta(X^{i_1}_0 t^{-1/2}, \ldots,X^{i_n}_0 t^{-1/2})  
\prod_{k=1}^{n} G_t(X^{i_k}_0,y_{n-k+1}) \right] \, \Delta(yt^{-1/2}) \, dy_1 \ldots dy_n \\
& \geq & C(n,L) t^{-n/2} \, P[\Omega_0(t)] \, \Delta(yt^{-1/2}) \, dy_1 \ldots dy_n
\end{eqnarray*}
where we have used (\ref{lbparticles}) and $|y_j| \leq L t^{1/2}$ in the final inequality.
We have also used $C(n,L,\ldots)$ to denote a finite non-zero quantity, depending only
on the quantities listed, but whose exact value is unimportant and may change from line to line.
It remains only to bound $P[\Omega_0(t)]$ from below, 
independently of $t \geq t_0$. This clearly holds  under the two sets of
assumptions described in section \ref{s1}, and in particular for Poissonian 
initial conditions.

For the upper bound in the asymptotic (\ref{eq:mainconj}) we estimate the probability that  
there exist particles at $dy_1,\ldots,dy_n$ at time $2t$ by conditioning on the set 
$\{X^i_t: i \in \N\}$ of positions of the particles at time $t$. For the desired 
particles to exist at time $2t$, one of the events
\[
\Omega_{i_1, \ldots,i_n}(t) = \left\{ \mbox{$X^{i_k}_{2t} \in dy_k$ and
the paths $(X^{i_k}_s:s \in [t,2t])$ do not intersect for $k=1,\ldots,n$} \right\},
\] 
for some $i_1<\ldots<i_n$, must occur.
Applying the Markov property at time $t$ and the upper bounds in Lemma
\ref{l1} we find, for all $t>0$,  
\begin{eqnarray*}
&& \hspace{-.4in} 
P \left[ \mbox{there exist particles at $dy_1,\ldots,dy_n$ at time $2t$} \right] \\
& \leq & \sum_{i_i<\ldots<i_n} P \left[\Omega_{i_1, \ldots,i_n}(t)\right] \\
& \leq & c_n \sum_{i_i<\ldots<i_n} E \left[ 
\Delta(X^{i_1}_t t^{-1/2}, \ldots,X^{i_n}_t t^{-1/2})
\prod_{k=1}^{n} G_t(X^{i_k}_t,y_k) \right] 
 \Delta(yt^{-1/2}) \, dy_1 \ldots dy_n \\
& \leq & C(n) \left( \int_{\R^n} 
G_t(x,y) \Delta(xt^{-1/2}) \rho_n(x;t) dx \right) \Delta(yt^{-1/2})
\, dy_1 \ldots dy_n \quad \mbox{see (\ref{usingrho})} \\
& \leq & C(n,L) \, t^{-n/2} \,  \Delta(yt^{-1/2}) \, dy_1 \ldots dy_n. 
\end{eqnarray*}
The final inequality uses the substitution $x \to xt^{1/2}$ and the bound
from Lemma \ref{l2}.
\section{Proof of Lemma \ref{l1}, the Karlin-Macgregor bounds} \label{s3}
First note that Brownian scaling implies that 
$G^{KM}_t(x,y) = t^{-n/2} \, G^{KM}_1(xt^{-1/2},yt^{-1/2})$,
and we need only prove the lemma for $t=1$.  
Factoring out common terms from the rows in (\ref{KM}) we find that
\[
G_1^{KM}(x,y) = (2 \pi)^{-n/2} \det (E(x,y))
\, \prod_{i=1}^{n} e^{-(x_i^2 + y_i^2)/2}
\]
where $E(x,y)$ is the matrix with entries $\exp(x_i y_j)$.
We now use 
the Harish-Chandra-Itzykson-Zuber formula \cite{Zub}, which is widely
used in random matrix theory. This states that 
\[
 \det (E(x,y)) = c_n \Delta(x) \Delta(y) 
\int_{\U(n)}  \exp \left( \tr(UX U^{\dagger}Y) \right) \mu(dU),
\]  
where $\tr(A)$ is the trace of a matrix $A$, 
$\mu(dU)$ is (normalized) Haar measure on the unitary group $\U(n)$
and $X,Y$ are the diagonal matrices with entries $x_1,\ldots,x_n$ and $y_1,\ldots,y_n$.
A short proof of this formula is found in \cite{Mehta} (appendix A.5). 

The fact that $U$ is unitary ensures that the function $F(U)= \tr(UX U^{\dagger}Y) $ is real 
valued (recall that $U^{\dagger}$ is the conjugate transpose of $U$ so that 
$U U^{\dagger}=U^{\dagger}U=I$).
By the compactness of $\U(n)$, $F(U)$ achieves its maximum and minimum values and we claim
that
\begin{equation} \label{claim1}
\min_{U \in \U(n)} F(U) = \sum_{i=1}^{n} x_{n-i+1} \, y_i, \quad
\max_{U \in \U(n)} F(U) = \sum_{i=1}^{n} x_{i} \, y_{i}.
\end{equation}
Bounding the integral in the Harish-Chandra-Itzykson-Zuber formula using these
maxima and minima in the integrand leads to the bounds (\ref{KMbounds}).

It remains only to prove the claim (\ref{claim1}), for which
we will find all the stationary points of the smooth function $F$ on the manifold $\U(n)$. 
A matrix $V$ is in the tangent space to $\U(n)$ at the point $U$ precisely when 
$UV^{\dagger} + VU^{\dagger} = U^{\dagger}V + V^{\dagger}U = 0$. 
The first derivative of $F$ in the direction of $V$ in the tangent space is given by  
\begin{eqnarray*}
D_V F(U) & = & \tr(UXV^{\dagger}Y) + \tr(VXU^{\dagger}Y) \\
& = & \tr(V^{\dagger}YUX) + \tr(VXU^{\dagger}Y) \quad 
\mbox{(since $\tr(AB) = \tr(BA)$ for Hermitian $B$)} \\
& = & \tr(UV^{\dagger}YUXU^{\dagger}) + \tr(VXU^{\dagger}Y) \quad 
\mbox{(since $\tr(UAU^{\dagger}) = \tr(A)$)} \\
& = & - \tr(VU^{\dagger}YUXU^{\dagger}) + \tr(VU^{\dagger}UXU^{\dagger}Y) 
\quad \mbox{(using $UV^{\dagger} = -  VU^{\dagger}$)} \\
& = & \tr \left( VU^{\dagger} \, [UXU^{\dagger},Y] \right) 
\end{eqnarray*}
(where $[A,B]= AB-BA$ is the commutator).
The map $V \to VU^{\dagger}$ bijectively maps the tangent space onto the 
Lie group $Lie(\U(n))$ defined by $\{W:W^{\dagger}+W=0\}$.
It is straightforward to check that $\tr(WA) =0$ for all $W \in Lie(\U(n))$ implies that
$A=0$. The vanishing of $D_V F(U)$ therefore implies that $[UXU^\dagger ,Y]=0$. 
This in turn implies, since $Y$ 
is diagonal with distinct eigenvalues, that $UXU^\dagger$ is a diagonal matrix. The matrix 
$UXU^\dagger$ has the same 
distinct eigenvalues as $X$, so that the entries on the diagonal of $UXU^{\dagger}$
must be a permutation of $x_1, \ldots,x_n$. Hence the value of $F(U)$ at a stationary point is
$\sum_{i=1}^{n} x_{\pi(i)} y_i$ for some permutation $\pi$. We must show that this value is
maximized when $\pi=\mbox{Id}$ the identity permutation and minimized at the involution
$\pi(i)=N+1-i$. Argue by contradiction and suppose that $\pi \to \sum_{i=1}^{n} x_{\pi(i)} y_i$ is 
maximized over permutations at some $\pi_0 \neq \mbox{Id}$. Then there exist $i<j$ for which
$\pi_0(i)>\pi_0(j)$. Let $\pi^{(jk)}$ be the transposition permutation, swapping the $j$th and $k$th 
elements. Then an increased value is found by taking the composition $\pi_1= \pi^{(jk)} \circ \pi_0$  
permutation, since
\[
\sum_{i=1}^{n} x_{\pi_1(i)} y_i - \sum_{i=1}^{n} x_{\pi_0(i)} y_i 
= (x_{\pi_0(i)} -x_{\pi_0(j)}) (y_j-y_i) > 0,
\]
using the ordering $x_1< \ldots <x_n$, $y_1< \ldots<y_n$. The minimum 
follows from a similar argument.
\section{Proof of Lemma \ref{l2}, the $n$-point density bound} \label{s4}
We start with a simple 
construction of coalescing particles as follows.
Fix $x_1<x_2< \ldots$. Take I.I.D. Brownian motions $(B^i:i \geq 1)$ and set
$\hat{X}^{x_i}_t = x_i+B^i_t$, for $t \geq 0$, the  
non-interacting Brownian paths.
Then define $X^{x_1}_t = \hat{X}^{x_1}_t$ and inductively for $i \geq 2$
\[
X^{x_i}_t = \left\{ \begin{array}{ll}
\hat{X}^{x_i}_t & \mbox{for $t \leq \tau(x_i,x_{i-1})$,} \\
X^{x_{i-1}}_t & \mbox{for $ t > \tau(x_i,x_{i-1})$,}
\end{array} \right.
\]
where $\tau(x_i,x_{i-1}) = \inf \{t: \hat{X}^{x_i}_t = X^{x_{i-1}}_t\} 
= \inf \{ t: X^{x_i}_t = X^{x_{i-1}}_t\}$. 
The path $t \to X^x_t$ gives the position of the particle that started at $x$.
This construction shows that the process
starting from infinitely many particles can be approximated by
the process starting with finitely many particles in an increasing way.
It will therefore be sufficient to bound $\rho_n$ over all finite starting positions. 
We write $P_{x_1,x_2,\ldots}$ to indicate the starting positions.

The fact that $\rho_1(y) \leq (\pi t)^{-1/2}$, for any
initial particle configuration, is well known. In section 9 of \cite{Evansetal} 
it states that the time reversal duality formula for coalescing 
Brownian motions for $a<b$ is given by
\[
 P_{x_1,\ldots,x_m} \left[ \{X^{x_i}_t: i \leq m\} \cap [a,b] \neq \emptyset \right] 
=
 P_{a,b} \left[ [X^a_t,X^b_t] \cap \{x_1,\ldots,x_m\} \neq \emptyset \right]. 
\]
We can take the supremum of the above duality over all starting 
configurations to obtain 
\[ 
\sup_{x_i} P_{x_1,\ldots,x_m} \left[ \{X^{x_i}_t: i \geq 1\} \cap [a,b] \neq \emptyset \right] 
= 
P_{a,b} \left[ X^a_t < X^b_t \right]  
\]
which can be explicitly calculated and is bounded by $(\pi t)^{-1/2} (b-a)$.
We give an elementary approach to this bound that avoids duality in Appendix \ref{appendixA}. 

For the extension to $\rho_n$ given in Lemma \ref{l2} 
we fix $a_1<b_1<a_2< \ldots <a_n < b_n$. Define
\[
\Omega_j = \left\{ 
\mbox{there exists a particle in $[a_j,b_j]$ at time $t$} 
\right\}.
\]
We shall show by induction that
\begin{equation} \label{l2-}
P_{x_1,\ldots,x_m} \left[ \Omega_1 \cap \ldots \cap \Omega_n \right]
\leq \prod_{i=1}^n \frac{(b_i-a_i)}{\pi t^{1/2}}.
\end{equation}
Decompose $\Omega_1$ as the disjoint union $\Omega_1 = \cup_k \Omega_{1,k}$ where
\[
\Omega_{1,k} = \left\{ 
\mbox{$X^{x_i}_t < a_1$ for $i=1,\ldots,k-1$ and $X^{x_k}_t \in [a_1,b_1]$} 
\right\}.
\]
Then conditional on $\Omega_{1,k}$ the processes $(X^{x_i}:i \geq k+1)$ evolve as coalescing 
Brownian motions with an extra lower absorbing boundary along the path
$t \to X^{x_k}_t$. The construction above shows that this extra absorbing 
path only lowers the probability that any of the paths 
$(X^{x_i}:i \geq k+1)$ reach the intervals $[a_2,b_2], \ldots, [a_n,b_n]$.
Thus, applying Bayes formula,
\begin{eqnarray*}
P_{x_1,\ldots,x_m} \left[ \Omega_1 \cap \ldots \cap \Omega_n \right] 
& = & \sum_{k=1}^m P_{x_1,\ldots,x_m} \left[ \Omega_2 \cap \ldots \cap \Omega_n \, | \, \Omega_{1,k} 
\right] P_{x_1,\ldots,x_m} \left[ \Omega_{1,k} \right] \\
& \leq & \sum_{k=1}^m P_{x_{k+1},\ldots,x_m} \left[ \Omega_2 \cap \ldots \cap \Omega_n \right] 
P_{x_1,\ldots,x_m} \left[ \Omega_{1,k} \right] \\
& \leq & \sum_{k=1}^m P_{x_{1},\ldots,x_m} \left[ \Omega_2 \cap \ldots \cap \Omega_n \right] 
P_{x_1,\ldots,x_m} \left[ \Omega_{1,k} \right] \\
&=& P_{x_{1},\ldots,x_m} \left[ \Omega_2 \cap \ldots \cap \Omega_n \right] 
P_{x_1,\ldots,x_m} \left[ \Omega_1 \right]. 
\end{eqnarray*}
and the result (\ref{l2-}), and hence Lemma \ref{l2}, follows. 
\begin{appendix}
\section{An elementary approach to the one point density}\label{appendixA}
This elementary approach may be of some interest for systems that do not 
have a duality relation.
Fix $a\in [0,1]$ and let $p(x_1,\ldots,x_n)$ be the probability that,
starting some coalescing Brownian motions $x_1,\ldots,x_n$, 
there is a particle at time $t$ in the interval $[0,a]$. 
We aim to show, by induction on $n$, that
\begin{equation} \label{inductivehypothesis}
p(x_1,\ldots,x_n) \leq C t^{-1/2} a
\end{equation}
for some constant $C<\infty$.
We suppose $x_1< \ldots < x_{n+1}$. Condition on the path $t \to X^{(n+1)}_t$
of the path started at $x_{n+1}$. We can consider three cases: if
$X^{(n+1)}_t < 0$ then we are sunk since no other path can ever enter
$[0,a]$ due to coalescence; if $X^{(n+1)}_t \in [0,a]$ then we are done;
finally if $X^{(n+1)}_t > a$ then we still need
some of the particles from $x_1, \ldots,x_n$ to enter $[0,a]$, but 
there is an absorbing boundary along the path of $X^{(n+1)}_t$. This absorbing boundary 
can only lower the chance of getting particles where we want. 
This splitting of possibilities leads to 
\begin{equation}
 p(x_1,\ldots,x_{n+1}) 
 \leq  C t^{-1/2} a P[X^{(n+1)} >a ]  + 
P[X^{(n+1)} \in [0,a]]. 
\end{equation}
It is straightforward to find  
$c_0< \infty$ so that for a Brownian motion $(X_t)$, and for $L \geq 1,t>0$,
\[
\frac{P[X_t \in [0,a] \, | \, X_0=x ]}{P[X_t < a  \, | \, X_0 = x]} \leq c_0 L a t^{-1/2}
\quad \mbox{whenever $|x| \leq L t^{-1/2}$.}
\] 
Suppose now that $x_i \in [-Lt^{1/2},Lt^{1/2}]$ for all $i$. Using this bound the induction
argument works with the choice $C=c_0 L$. In particular 
we have shown that
\[
\sup_{|x_i| \leq L t^{1/2}}  p(x_1,\ldots,x_m)
 \leq  c_0 L t^{-1/2} a. 
\]
Applying this with the choice $L=O(\sqrt{\log(t)})$, and using a simple bound
that any particle starting outside $[-Lt^{1/2},Lt^{1/2}]$ can reach the 
intervals $[0,a]$ (we need only control particles starting exactly at
$\pm Lt^{1/2}$ by coalescence) we find that
\[
\sup_{x_i}  p(x_1,\ldots,x_m) 
 \leq C t^{-1/2} (\log(t))^{1/2} a
\]
where the supremum is over all possible initial configurations. 
One can remove the unwanted logarithm term by a blocking argument, but
this is already longer than the elegant duality argument, which has the extra advantage of 
achieving the optimal constant.
\section{A note on the existence of $\rho_n$}
We briefly consider for which measure $\rho_n(x;t)$ acts as a density. 
List the (disjoint) positions of the particles at time $t$ as $\{X^i_t:i \geq 1\}$. 
Define a random measure $\mu^n_t$ on $\R^n$ by
\[
\mu^n_t = \sum_{i_1} \sum_{i_2 \neq i_1} \cdots \sum_{i_n \neq i_1,\ldots,i_{n-1}}
\delta (X^{i_1}_t,\ldots,X^{i_n}_t)
\]
where $\delta(x)$ is a point mass at $x$.  
Let $C_n(\delta)$ denote the cube $[0,\delta)^n$. The cubes $k + C_n(2^{-m})$, for $k \in (2^{-m}\Z)^n$, 
partition $\R^n$. Take $A \subseteq \R^n$ open. For an initial condition with only finitely many particles, 
the measure $\mu^n_t$ has only finitely many atoms which do not lie on $\cup_{i \neq j} \{x_i=x_j\}$.
Therefore we have the increasing limit 
\[
\lim_{m \to \infty} \sum_{k \in (2^{-m}\Z)^n} \I (\mu^n_t(k+C_n(2^{-m})) >0)
= \mu^n_t(A),  
\]
where the sum over $k$ is restricted to those terms for which $k_1,\ldots,k_n$ are distinct.
The event $\mu^n_t(k+C_n(2^{-m}))>0$ can be rewritten as
$[k_j,k_j+2^{-m}] \cap \{X^i_t:i \geq 1\} \neq \emptyset$ for $j=1,\ldots,n$. 
Taking expectations and applying (\ref{l2-}) we find
that $ E[ f d\mu^n_t ] \leq (2\pi t)^{-n/2} \int f dx $
when $f = \I(A)$. By an approximation argument this also holds for all
measurable $f$ and by monotonicity it holds for any initial distribution of particles. 
We may then define $\rho_n(x;t) = \rho(x_1, \ldots,x_n;t)$ as the density of the 
measure $E[\mu^n_t]$, that is the Radon Nicodym derivative with respect to Lebesgue 
measure, so that
\begin{equation} \label{usingrho}
\int_{\R^n} f(x) \rho_n(x;t) dx = E \left[  
\sum_{i_1} \sum_{i_2 \neq i_1} \cdots \sum_{i_n \neq i_1,\ldots,i_{n-1}}
f(X^{i_1}_t, \ldots, X^{i_n}_t) \right].
\end{equation}
The Lebesgue differentiation theorem implies that 
$ \rho(x;t) = \lim_{\delta \to 0} \delta^{-n} E [ \mu^n_t(x+C_n(\delta)]$
for almost all $x$. 
Choosing $x_1,\ldots,x_n$ disjoint we obtain 
\[ 
\rho(x_1,\ldots,x_n;t) = \lim_{\delta \to 0} \delta^{-n} E \left[ N_t([x_1,x_1+\delta]) \ldots
N_t([x_n, x_n + \delta]) \right]
\]
where $N_t(A)$ is the number of particles inside $A$ at time $t$. However we claim we may replace 
$N_t([x_1, x_1+\delta])$ by the indicator $\I( \mbox{$[x_1, x_1+\delta]$ contains a particle})$ 
in this limit. The error when doing this occurs if there are two or 
more particles in $[x_1,x_1+\delta]$ and is dominated by 
$ E[ \mu^{n+1}_t( (x_1,x_1,x_2,\ldots,x_n)+C_{n+1}(\delta))]$. 
This error is therefore of order $\delta^{n+1}$, so that 
replacing $N_t([x_i, x_i+\delta])$ in this way for each $i=1,\ldots,n$ 
we reach the usual definition of $\rho_n(x;t)$ as, for distinct $x_1,\ldots,x_n$,
\[
\rho(x_1, \ldots,x_n;t) = \lim_{\delta \to 0} \delta^{-n} P\left[ \mbox{all the intervals
$[x_i,x_i+\delta]$ for $i=1,\ldots,n$ are non-empty}\right].
\] 
\end{appendix}

\end{document}